\documentclass[smallextended,envcountsect]{svjour3} 
\smartqed
\usepackage{graphicx}
\usepackage{amsfonts}
\usepackage{amsmath}
\usepackage[usenames,dvipsnames]{pstricks}
\usepackage{epsfig}
\usepackage{bbding}

\newcommand{\bc}{\begin{center}}
\newcommand{\br}{\begin{right}}
\newcommand{\ec}{\end{center}}
\newcommand{\be}{\begin{equation}}
\newcommand{\ee}{\end{equation}}

\newcommand{\vl}{\mid}

\newcommand{\rar}{\rightarrow}

\newcommand{\p}{\partial}

\newcommand{\bmax}{\mbox{max} \,}

\newcommand{\eps}{\varepsilon}

\newtheorem{thm}{Theorem}[section]

\newtheorem{defi}{Definition}[section]

\newtheorem{rem}{Remark}[section]

\journalname{COMMUNICATIONS IN OPTIMIZATION THEORY}

\begin{document}

\title{Nash equilibrium points and their finding for nonsmooth case}


\titlerunning{Nash equilibrium points and their finding }

\author{Igor M. Prudnikov}

\institute{Igor Mihailovich Prudnikov \at
             Scientific Center of Smolensk Federal Medical  University, Smolensk, Russia, 214000\\
              \email{pim\underline{ }10@hotmail.com}
}

\date{Received: date / Accepted: date}

\maketitle

\begin{abstract}

The purpose of this paper is to develop a numerical method for
finding an equilibrium point in a model, in which the loss
function of each object (subject) is described by a convex
function with respect to one of its variables. Such models are
found in medicine, economics, game theory, and biology.

For the more complex case, with nonsmooth functions describing the
state of each element of the system as damage, loss, or gain, the
Steklov average integrals are used that turn nonsmooth functions
into smooth ones.

Numerical methods for finding equilibrium points in the more
general non-smooth case are constructed. In the process of
optimization, the diameters of the sets, over which the averaging
takes place, are decreased in accordance with the optimization
steps.

All limit points are proved to be equilibrium points. Under some
conditions, the convergence rate can be estimated using the
Kantorovich theorem. The necessity to develop new methods for
finding Nash equilibrium points in the nonsmooth case is
concluded.

\end{abstract}



\section{Introduction}

Let the physical or economic state of a system be described by $m$
loss functions $f_1(x_1, x_2, \dots, x_m):\mathbb{R}^m \rar
\mathbb{R}$, $f_2(x_1, x_2, \dots, x_m): \mathbb{R}^m \rar
\mathbb{R}$, $\dots$, $f_m(x_1, x_2, \dots, x_m):\mathbb{R}^m \rar
\mathbb{R}$ depending on $m$ variables $x_1, x_2, \dots, x_m,$
where $\mathbb{R}^m$ is $m$-dimensional Euclidean space. Then an
equilibrium point is a state $x^*_1, x^*_2, \dots, x^*_m$ for
which changing any $x^*_j$ leads to an increase in the
corresponding  function $f_j(\cdot)$, i.e.  \be f_j(x^*_1, x^*_2,
\dots,x^*_j,\dots, x^*_m) \leq f_j(x^*_1, x^*_2, \dots, x^*_{j-1},
x_j, x^*_{j+1}, \dots, x^*_m) \label{equilibpoints1} \ee

Equilibrium states were introduced into economics by J. Nash. In
1950-1953, his articles proving the existence of equilibrium
points were published \cite{Nash1}-\cite{Nash4}.

The problem of equilibrium point finding in biology or economics
is closely related to game theory and is of practical importance.
Equilibrium points arise from interspecific competition in biology
and intercompany competition in economics. Equilibrium points in
medicine are homeostasis points \cite{gomeostasis}. These are
points of balance between various states of the human body, e.g.
blood pressure, temperature, blood cholesterol level, pulse rate.
Some balance is achieved between different pills when we take
medication.

Here we consider non-cooperative games of $m$ players, none of
whom can influence other players' behavior (strategies). A player
$i$ chooses independently a pure strategy $x_i$ from a compact
convex set $S_i$, such that he minimizes his  loss  function
$f_i(\cdot)$.

Consider a vector $x=(x_1, x_2, \dots, x_m) \in \mathbb{R}^m,$
called a multistrategy  and comprised of the pure strategies $x_i
\in S_i$. We assume that the vector $x$ belongs to the compact
convex set $S=S_1 \times S_2 \times \dots S_m \in \mathbb{R}^m,$
which is the Cartesian product of the compact sets $S_i, i \in
1:m$, and $\mbox{int } S \neq \emptyset$.

\begin{defi} A multistrategy $x^* = ( x^*_1, x^*_2, \dots,
x^*_m) \in  S $ of a non-cooperative game is called a
non-cooperative equilibrium if the inequality
(\ref{equilibpoints1}) is true for every $j \in 1:m$ and $x_j \in
S_j$.
\end{defi}

The equilibrium point definitions in game theory, medicine, and
economics are similar. In 1950, J. Nash proved the following
theorem.

\begin{thm} \cite{Nash1}. Let  $S_i$ be a compact convex  set for
any $i \in 1:m$ and $f_i(\cdot)$ be convex  with respect to $x_i
\in S_i$. Then there is a non-cooperative equilibrium in a
non-cooperative game with $m$ players.
\end{thm}

The aim of the paper is to develop numerical methods for finding
equilibrium points in the nonsmooth case. While, to our knowledge,
papers describing numerical methods for finding equilibrium points
in special cases \cite{minarchenko} have been published, no papers
describing numerical methods for finding equilibrium points in the
general case are found.

\vspace {0.5cm}

\section{\bf Discussion of the problem}

\vspace {0.5cm}

We describe a method of searching for an equilibrium state,
provided that the functions $ f_i (x _ {- i}, x_i): S \rar R $,
where $ x_{- i} = (x_1, x_2, \dots, x_ {i-1}, x_ {i + 1}, \dots,
x_m) $, are convex with respect to  $ x_i $. We assume that  the
inclusion $ \{x \in \mathbb{R}^m \vl f_i (x) <f_i (x_0) \} \subset
\mbox {int } S $ is true for any $ i \in 1: m $, in which $ x_0 $
is a starting point. Here $\mbox {int} S $ means the interior of
the set $ S $.

Let us denote the coordinate vectors  by $ e_1 = (1,0, 0, \dots,
0), e_2 = (0,1,0, \dots, 0), \\ \dots, e_m = (0,0,0, \dots, 1)$.
It is known that coordinate descent method admits no convergence
for nonsmooth functions \cite{demminmax}. Therefore, we will use
the ideas from \cite{proudintegapp1}.

\begin{center}
{\bf Equilibrium points finding algorithms for the smooth case}
\end{center}

One way is to use gradient and second-order methods for the smooth
case, i.e., when $f_i(\cdot), i \in 1:m, $ are differentiable
functions with respect to the variables $x_j, j \in 1:m. $ Denote
the partial derivative of the function $f_i(\cdot)$ with respect
to the variable $x_i$ by $f'_{i,x_i}(\cdot)$.

Consider the vector function $\Theta(\cdot): \mathbb{R}^m \rar
\mathbb{R}^m :$ \be \Theta(x)=
\begin{pmatrix}f'_{1,x_1}(x)\\f'_{2,x_2}(x)\\
\ldots\\f'_{m,x_m}(x)
\end{pmatrix}
\label{equilibpoints01} \ee We search a vector $x^*=(x^*_1, x^*_2,
\dots, x^*_m), $ for which $\Theta(x^*)=0.$ It is clear that  the
vector $x^*$ is an equilibrium point.

From expansion accurate to the higher order terms
$$
\Theta(x+\Delta x)=\Theta(x)+ \Theta'(x) \Delta x +o(\Delta x),
$$
in which
$$
\lim_{\Delta x \rar 0} \frac{o(\Delta x)}{ \| \Delta x \|} =0,
$$
we obtain the value for the step $\Delta x$. Suppose
$$
\Theta(x+\Delta x) \approx \Theta(x)+ \Theta'(x) \Delta x = 0.
$$
From here we obtain
$$
\Delta x = - (\Theta'(x))^{-1} \Theta(x),
$$
where \be \Theta'(x) = \begin{pmatrix}
f''_{1,x_1,x_1} & f''_{1,x_1,x_2} & \cdots & f''_{1,x_1,x_m} \\
f''_{2,x_2,x_1} & f''_{2,x_2,x_2} & \cdots & f''_{2,x_2,x_m} \\
\vdots & \vdots & \ddots & \vdots \\
f''_{m,x_m,x_1} & f''_{m,x_m,x_2} & \cdots & f''_{m,x_m,x_m} \\
\end{pmatrix}
\label{equilibpoints02} \ee

\vspace{0.3cm}
  {\bf Algorithm 1 (The Newton's method for twice continuously
differentiable  functions  $f i(\cdot)$)} \vspace{0.3cm}

At each step $k$ we find \be \Delta x_k=-
(\Theta'(x_k))^{-1}\Theta(x_k). \label{equilibpoints2} \ee We set
$x_{k+1}=x_k+2^{-l} \, \Delta x_k,$ in which $l$ is the smallest
number from the set $M=\{ 0,1,2, \cdots \}$ for which the
inequality $\| \Theta(x_{k+1}) \| < \| \Theta(x_k) \|$ is correct.
It is easy to prove that there exists a number $l$ for which the
inequality $\| \Theta(x_{k+1}) \| < \| \Theta(x_k) \|$ is correct.
Repeat the process as long as $\| \Theta(x_k) \| \leq \varepsilon,
$ in which $ \varepsilon $ is a positive small number.

Let  the inequality \be L_1 \| \Delta x \| \leq\| \Theta'(x)
\Delta x \| \leq L_2 \| \Delta x \| \label{equilibpoints12} \ee
hold true for some $L_1, L_2 >0$ and any $\Delta x$. We assume
that $x$ belongs to a small neighborhood of an equilibrium point
$x^*$, in which $\Theta(x^*)=0$ and the optimization process takes
place with full step $\Delta x_k$ i.e. $l=0$. Then it is possible
to obtain an estimation of the convergence rate of the Newton's
method. We have
$$
\| \Theta (x_{k+1}) \| = \| \Theta (x_k+\Delta x_k)-\Theta(x^*)
\|= \| \Theta'(\xi) \Delta x_k \|,
$$
in which $\xi$ is a point on the line, connecting $x_{k+1}$ and $
x^*$, and between them. Let's substitute the expression for
$\Delta x_k$ from (\ref{equilibpoints2}). Therefore, we obtain
$$
\| \Theta'(\xi) \Delta x_k \| = \| \Theta'(\xi)
(\Theta'(x_k))^{-1} \Theta(x_k) \|.
$$
Due to the continuity of the matrix $\Theta'(\cdot)$ and the fact
that at each step $k$ according to the choice of the step
$\Theta(x_{k})=o(\Delta x_{k-1})$, and also considering the
assumption (\ref{equilibpoints12}) we obtain a chain of
inequalities
$$
L_1 \| \Delta x_k \|   \leq  \| \Theta'(\xi) \Delta x_k \| =
\|\Theta'(\xi) (\Theta'(x_k))^{-1} \Theta(x_k) \| \leq c_k
\frac{\| \Delta x_{k-1} \|}{N_{k-1}}.
$$
Here $c_k=\| \Theta'(\xi) ( \Theta'(x_k))^{-1} \|$. We have
$$
\| \Theta(x_k) \| = \| o(\Delta x_{k-1}) \| \leq \frac{\| \Delta
x_{k-1} \|}{N_{k-1}}
$$
in which $N_{k-1} = N(\Delta x_{k-1}) \rar_k \infty$ as $ \|
\Delta x_{k-1} \| \rar 0 $. The latter follows from the definition
of the  infinitesimal function $o(\cdot)$. We have
$$
\| o(\Delta x) \| \leq \eps(\Delta x) \| \Delta x \|,
$$
in which $ \eps(\Delta x) \rar 0 $ as $\Delta x \rar 0. $
Therefore, we can put
$$
N_{k} = N(\Delta x_{k})= \frac{1}{\eps(\Delta x_{k})}.
$$
From the inequalities written we get
$$
\| \Delta x_k \| \leq c_k \frac{\|\Delta x_{k-1} \|}{L_1 \,
N_{k-1}} = q_k \| \Delta x_{k-1} \|.
$$
Superlinear convergence follows from here, since
$$q_k= \frac{c_k}{L_1 \, N_{k-1}} \rar_k 0.$$ We finally obtain
\begin{thm}
Let the assumption (\ref{equilibpoints12}) hold true for the twice
continuously differentiable functions $f_i(\cdot), i \in 1:m.$
Then the Newton's method will converge with superlinear velocity
in a small neighborhood of $ x^*$. \label{thmequilibpoints0}
\end{thm}
This optimization process requires the existence of continuous
second mixed derivatives with respect to the variables $x_i, x_j,
i,j \in 1:m,$ of the functions $f_i(\cdot), i \in 1:m,$ and the
existence of the inverse matrix $(\Theta'(x_k))^{-1}$ at any step
$k$. Unfortunately, the assumption (\ref{equilibpoints12}) does
not always hold. Moreover, the theorem {\ref{thmequilibpoints0} is
true in a small neighborhood of the point $x^*$ which is to be
reached.

\vspace{0.5cm}

\section{\bf Solution of the problem for the nonsmooth case}

\vspace{0.5cm}

Let us use the ideas of the paper \cite {proudintegapp1}. We
assume, that $ f_i (\cdot), i \in 1:m, $ are Lipschitz functions
with constants $L_i$ i.e.
$$
\| f_i(u) - f_i(v) \| \leq L_i \| u - v \|
$$
for all $u,v \in \mathbb{R}^m$. We construct functions \be
\varphi_i (x) = \frac {1} {\mu (D)} \int_D f_i (z + x) d z. \label
{equilibpoints11} \ee in which $ D $ is an arbitrary convex
compact set, $ 0 \in \mbox {int} D, \mu (D)> 0 $ is the Lebesgue
measure of the set $ D$ and the integral is the Lebesgue integral.
It is not difficult to verify that the function $ \varphi_i
(\cdot) $ is convex with respect to the variable $ x_i.$

The function $f_i(\cdot)$ has the partial derivative with respect
to $x_i$ almost everywhere (a.e.) on the set $S$. In
\cite{proudintegapp1} it was proven that the function
$\varphi_i(\cdot) $ is continuously differentiable with respect to
the variable $ x_i $. The partial derivative $\varphi_i(\cdot) $
with respect to $ x_i $ can be calculated by the formula
\cite{proudintegapp1}
$$
\varphi'_{i,x_i}(x)=\frac{\p \varphi_i(x)}{\p
x_i}=\frac{1}{\mu(D)} \int_D \frac{\p f_i (x+y)}{\p x_i} d y.
$$
The functions $\varphi_i(\cdot), i \in 1:m,$ have an equilibrium
point according to Nash's theorem.

Substitute  $f_i(\cdot), i \in 1:m,$ for $\Phi_i(\cdot):
\mathbb{R}^m \rar \mathbb{R}, i \in 1:m,$ defined by
$$
\Phi_i(x)= \frac{1}{\mu(D)} \int_D \varphi_i(x+y)d y,
$$
in which the functions $\varphi_i(\cdot), i \in 1:m,$ and the set
$D$ are defined above (see (\ref{equilibpoints11})). We take the
integral of the integral, since in this way we obtain the twice
continuously differentiable functions $\Phi_i(\cdot), i \in 1:m,$
and the stationary points of $\Phi_i(\cdot)$ are $\eps(D)-$
stationary points of $f(\cdot)$ \cite{proudintegapp1}.

Since $\varphi_i(\cdot)$ is a Lipschitz \cite{proudintegapp1}, we
will have \be \Phi_i'(x)= \frac{1}{\mu(D)} \int_D \varphi_i'(z+x
)d z. \label{intapp4} \ee.

We have proved that the functions $\Phi_i(\cdot), i \in 1:m,$ have
Lipschitz second derivatives \cite{proudintegapp1}. If $D$ is a
ball or a cube centered at zero with the diameter $d(D)$, then the
functions $\Phi_i(\cdot), i \in 1:m,$ have Lipschitz second
derivatives $\Phi''_i(\cdot)$ with constant \cite{proudintegapp1}
$$
L_i'=\frac{2 L_i}{d^2(D)}.
$$
We can apply the Newton's method to the functions $\Phi_i(\cdot),
i \in 1:m,$ to find the equilibrium points. In the process of
optimization we will consistently decrease the step $\lambda_k$
and the diameter $d(D_k)$ so that the inequality  \be
 \frac{\lambda_k }{d^2(D_k)}< \varepsilon_k,
\label{equilibpoints14} \ee was true for some sequence $\{
\varepsilon_k  \} $  in which $ \varepsilon_k \rar +0 $ as $k \rar
\infty. $  We will prove that the inequality
(\ref{equilibpoints14}) guarantees that any limit point of a
sequence obtained by  the Newton's method using the functions
$\Phi'_i(\cdot), \Phi''_i(\cdot), i \in 1:m,$  is an equilibrium
point of the functions $f_i(\cdot), i \in 1:m$.

\vspace{0.3cm}
 {\bf The Newton's method for finding equilibrium points for $f_i(\cdot)$,
 $i \in 1:m,$ using the functions $\Phi_i(\cdot)$ }

\vspace{0.3cm}

Calculate $\Theta(z_k)$ and $\Theta'(z_k)$
$$ \Theta(x)=
\begin{pmatrix}\Phi'_{1,x_1}(x)\\\Phi'_{2,x_2}(x)\\\ldots\\\Phi'_{m,x_m}(x),
\end{pmatrix}
$$
$$ \Theta'(x) = \begin{pmatrix}
\Phi''_{1,x_1,x_1} & \Phi''_{1,x_1,x_2} & \cdots &
\Phi''_{1,x_1,x_m} \\\Phi''_{2,x_2,x_1} & \Phi''_{2,x_2,x_2} &
\cdots & \Phi''_{2,x_2,x_m} \\\vdots & \vdots & \ddots & \vdots
\\\Phi''_{m,x_m,x_1} & \Phi''_{m,x_m,x_2} & \cdots &
\Phi''_{m,x_m,x_m} \\\end{pmatrix}
$$
accordingly to (\ref{intapp4}) for the twice differentiable
functions $\Phi_i(\cdot), i \in 1:m,$ when $D$ is a ball or a
cube.

Take a sequence of sets $\{ D_s \}, s=1,2,\dots $ with non-empty
interior, the diameters $d(D_s)$  of which tend to zero in $s \rar
\infty$. Let $D_s=B^m_{r_s}(0)=\{ v \in \mathbb{R}^n \vl \| v \|
\| \leq r_s \}$ for $r_s \rar +0$ and $s \rar \infty$. Let us
introduce for $i \in 1:m$ the following sequence of functions
$$
\varphi_{i,s}(x)= \frac{1}{\mu(D_s)} \int_{D_s} f_i(x+y)d y
$$
and \be \Phi_{i,s}(x)= \frac{1}{\mu(D_s)} \int_{D_s} \varphi_{i,s}
(x+y)d y \label{equilibpoints14b}. \ee  The difference between
(\ref{intapp4}) and (\ref{equilibpoints14b}) is that
(\ref{intapp4}) is written for a constant $D$, while
(\ref{equilibpoints14b}) is written for a set $D_s $ depending on
the parameter $s$.

Construct the functions $\Theta_s (\cdot)$ for the functions
$\Phi_{i,s}(\cdot), i \in 1:m, $ as written above. Let the
inequality $ \| \Phi''_{i,s} (\cdot)\| \leq {L}_{s}$ hold true in
which $ \Phi''_{i,s} (\cdot)$ is the matrix of the second mixed
derivatives. In \cite{proudintegapp1} it was proved that
$L_s=\frac{L}{d(D_s)},$ in which $L=\bmax_{i \in 1:m} L_i$.

It follows from here that, depending on the selected metric of the
space $ \mathbb {R} ^ m $, the norm $ \| \Theta'_s (\cdot) \| $ is
proportional to $ L_s $. Suppose $ \| \Theta'_s (\cdot) \| \leq
L_s. $

Define the vector-function $\tilde{ \Theta}_s(\cdot): \mathbb{R}^m
\rar \mathbb{R}^m $ as a function of $y$: \be \tilde{
\Theta}_s(y,x)= \Theta_s (y) +2 L_s ( y - x ).
\label{equilibpoints14a} \ee

Then we have the inequality for the matrix $\tilde{
\Theta}'_s(\cdot)$   \be {L_s} \| z \|^2 \leq (\tilde{
\Theta}'_s(x,x)z,z) \leq 3{L}_{s} \| z \|^2 \,\,\,\, \forall z \in
\mathbb{R}^n. \label{equilibpoints15} \ee

Let us construct the Newton's method for finding the roots of the
equation $\Theta_s (x) = 0 $ using the function $\tilde {\Theta}_s
(\cdot) $. We will use the rule of consistent reduction of the
length $\lambda_k $ of $k^{th}$ step and the diameter $d(D_k)$.

{\bc \bf Description of the Newton's method for finding for the
equilibrium points  using $\Phi_{i,s}(\cdot)$. \ec}


Let a point $x_k$ were constructed at the step $k$. Construct the
point $x_{k+1}$. Take by definition $\tilde{ \Theta}_{s,k}(\cdot)
= \tilde{ \Theta}_{s}(\cdot, x_k)$. The dependence of $s$ on $k$
will be written as $s=s(k)$.

We calculate $\Delta x_k=- (\tilde{\Theta}_{s,k}'(x_k))^{-1}
\tilde{\Theta}_{s,k}(x_k)$ at each step $k$. We set $x_{k+1} =
x_k+2^{-l} \, \Delta x_k,$ in which $l$ is the smallest number
from the set $M= \{ 0, 1, 2, \cdots \} $ for which $\|
\tilde{\Theta}_{s,k}(x_{k+1}) \| < \| \tilde{\Theta}_{s,k}(x_k)
\|$.

It is possible to prove that $\| \tilde{\Theta}_{s,k}(x_{k+1}) \|
< \| \tilde{\Theta}_{s,k}(x_k) \|$ for small $ \| \Delta x_k \|$
and $ \| \Delta x_k \| \rar 0$ as $k \rar \infty $ for fixed $s$
in a small neighborhood of the equilibrium point $x^*$. We assume
that we reach a small surrounding of the equilibrium point $x^*$
for big $k$ in which the process takes place with the full step
$\Delta x_k $.

If the inequality \be \frac{\| \Delta x_k \|}{d^2(D_{s(k)})} <
\eps_k \label{equilibpoints18a} \ee is fulfilled for a sequence
$\{ \eps_k \}, \eps_k \rar +0$ as $k \rar \infty$, then we
decrease the diameter  $d(D_{s(k)})$ of $D_{s(k)}$ and increase
$k, s=s(k)$.

The inequality (\ref{equilibpoints15}) holds true for $
\tilde{\Theta}'_{s(k),k}(\cdot) $ and all $s, k$. Firstly, we
prove that \be \lim_{k \rar \infty}{\Theta}_{s(k)}(x_k) = 0
\label{equilibpoints16} \ee and the sequence $\{ x_k \}$  has a
limit point $x^*$.

We have the expansion  of the function
$\tilde{\Theta}_{s(k),k}(\cdot) $ in the neighborhood of the point
$x_k$
$$ \tilde{\Theta}_{s(k),k}(x_{k+1}) = \tilde{\Theta}_{s(k),k}(x_k))
+ \tilde{\Theta}_{s(k),k}'(x_k)) \Delta x_k + o_{s(k),k} (\Delta
x_k).
$$
After substitution $\Delta x_k=-
(\tilde{\Theta}_{s(k),k}'(x_k))^{-1} \tilde{\Theta}_{s,k}(x_k)$ in
this expansion we obtain \be \tilde{\Theta}_{s(k),k}(x_{k+1})=
o_{s(k),k} (\Delta x_k). \label{equilibpoints17} \ee

Let us prove that $o_{s(k),k} (\Delta x_k)$ is an infinitesimal
function with respect to $\Delta x_k$ as $k \rar \infty$. Since
$\tilde{\Theta}_{s(k),k}(\cdot)$ was obtained from $
{\Theta}_{s}(\cdot) $ through adding the linear function,
$o_{s(k),k} (\cdot)$ is the same infinitesimal function in the
expansion of $ {\Theta}_{s}(\cdot) $ in the surrounding of $x_k$.
Now we will obtain the upper bound for $o_{s(k),k} (\cdot).$

The following expansion takes place
$$
{\Theta}_{s}(x_{k+1}) = {\Theta}_{s}(x_k) + {\Theta}'_{s}(x_k)
\Delta x_k + o_{s,k} (\Delta x_k).
$$
Since the function ${\Theta}_{s}(\cdot)$ is continuously
differentiable for each $s$, according to the midpoint theorem, we
have
$$
{\Theta}_{s}(x_{k+1}) - {\Theta}_{s}(x_k)= {\Theta}'_{s}(\xi)
(x_{k+1} - x_k) = {\Theta}'_{s}(\xi) \Delta x_k,
$$
in which $ \xi \in [x_k, x_{k+1}]. $ Let us substitute this
difference in the Taylor series and use the  Lipschitzness of
${\Theta}'_{s}(\cdot)$ with the constant $\frac{2 L}{d^2(D_s)} $.
Therefore, since $\Theta'_{s}(\cdot) $  is Lipschitz with the
constant $ \frac{2 L}{d^2(D_s)} $ \cite{proudintegapp1}, we obtain
$$
\| o_{s,k} (\Delta x_k) \| \leq \| ( {\Theta}'_{s}(\xi) -
{\Theta}'_{s}(x_k) ) \Delta x_k \| \leq \frac{2 L\| \Delta x_k
\|}{d^2(D_s)} \| \Delta x_k \|.
$$
From here \be \frac{\| o_{s,k} (\Delta x_k) \|}{\| \Delta x_k \|}
\leq \frac{2 L \| \Delta x_k \|}{d^2(D_s)}.
\label{equilibpoints18} \ee Hence, if \be \lim_{k \rar \infty}
\frac{\| \Delta x_k \|}{d^2(D_{s(k)})} = 0 \label{equilibpoints19}
\ee holds true during optimization, then uniform infinitesimality
of $o_{s(k),k} (\cdot)$ with respect to $s=s(k)$ and $k$  follows
from here. However,  we organize  our process in such a way that
the limit equality (\ref{equilibpoints19}) was correct.

The limit equality
$$
\lim_{k \rar \infty} L_{s(k)} \| x_{k+1} - x_k \| = \lim_{k \rar
\infty} L_{s(k)} \Delta x_k = \lim_{k \rar \infty}
\frac{L}{d(D_{s(k)})} \Delta x_k \leq
$$
\be \leq  \lim_{k \rar \infty} \frac{L}{d^2(D_{s(k)})} \Delta x_k
= 0. \label{equilibpoints19a} \ee follows from the inequality
(\ref{equilibpoints18a}) as we decrease $ d(D_{s(k)}) $ in the
process of optimization. The equality (\ref{equilibpoints16})
follows from (\ref{equilibpoints14a}), (\ref{equilibpoints17}) and
(\ref{equilibpoints19a}).

It follows from the upper semicontinuity of the Clarke
subdifferential and from the equality (\ref{equilibpoints16}) that
the sequence $\{ x_k \}$ converges to a limit point $x^*$, in
which $0 \in \p_{x_i} f_i(x^*)$ for all $i \in 1:m, $ i.e. $x^*$
is the equilibrium point.

All of the above stated is true if we reach a small neighborhood
of the equilibrium point. In order to do this, we are to use the
coordinate descent method  with some modifications for the
functions $\varphi_i(\cdot), i \in 1:m.$

Thus, the following theorem is proved.
\begin{thm}
Any limit points of the sequence obtained by Newton's method with
starting points from small neighborhoods of the equilibrium
points, are the equilibrium points if the equality
(\ref{equilibpoints19}) is satisfied in the process of
optimization for the convex with respect to $x_i$, Lipschitzian
functions $f_i(\cdot), i \in 1:m,$ respectively.
\label{thmequilibpoints3}
\end{thm}

The given Newton's method is also called the modified Newton's
method. It is possible to show that there is a majorant
Kantorovich function for any step $k$ \cite{kantorovichakilov}.
The step length and the convergence rate of the optimization
method are estimated under the conditions of consistency
(\ref{equilibpoints19}) and some conditions indicated in the
theorem given below. We can state  the convergence of the whole
sequence $\{ x_k \}$ under the below given conditions in the
theorem \ref{thmequilibpoints3}.

This is true for the reason that $\| \Delta x_k \| $ is compared
with the step length of the majorant function. The conditions of
the Kantorovich theorem \cite{kantorovichakilov}, pp. 689-690, are
fulfilled if we satisfy some requirements.

We will construct a sequence $\{ x_k \}$  converging to the
solution of the equation $\Theta(x^*)=0 $ for a ball $B^m_r(x_0)=
\{ y \in \mathbb{R}^m \vl \| y - x_0 \| \leq r \} $. Suppose
$Q_{s,0} = [\tilde{\Theta}'_{s,0}(x_0)]^{-1}$, $\| Q_{s,0}
{\Theta}_{s} (x_{0}) \| \leq A_{s}$ and $ \| Q_{s,0}
\tilde{\Theta}''_{s,0} (x) \| \leq B_{s} $ for any $x \in
B^m_r(x_0)$, $ \triangle x_k = -[\tilde{\Theta}'_{s,0}(x_k)]^{-1}
\tilde{\Theta}_{s,k}(x_k)$. We set $x_{k+1}=x_k+\triangle x_k$.
During the optimization process, we change $s=s(k)$ and  decrease
the diameter $d(D_s)$ so that the requirements of Theorem
\ref{thmequilibpoints4} were satisfied.

\begin{thm} We make the following assumptions: \\
There exists a linear operator $ Q_{s,0} = [\tilde {\Theta} '_ {s,
0} (x_0)] ^ {- 1} $ for
$ s = s (0) $ and $ k = 0 $. \\
If
$$
q_{s} = A_{s} B_{s} \leq q <\frac{1} {2}, \, \,
$$
is true for any $ s $ and the consistency condition
$$
\lim_{k \rar \infty} \frac{\Delta x_k}{d^2(D_{s(k)})}  = 0.
$$
is satisfied, then the equation $ {\Theta} (x) = 0 $ has a
solution $ x ^ * $ to which the Newton's method converges with the
rate \be \| x^* -x_{k} \| \leq \frac{1}{2^k} [2 q]^{2^k} C
\label{equilibpoints19b} \ee for a constant $C$.

The convergence rate of the modified Newton's method (for
$q<\frac{1}{2}$) is estimated by the following inequality \be \||
x ^ * - x_k \| \leq C (1- \sqrt {1-2q}) ^ {k + 1}, \, \,  k =
0,1,2, \dots \label{equilibpoints19d}   \ee
\label{thmequilibpoints4}
\end{thm}

\begin{rem}
The convergence rate proof follows with some changes from
\cite{kantorovichakilov}, p. 690, since the convergence rate
depends on the values $ q_k $ and $ A_k $. The first value is
limited by the value $ \frac{1}{2}$. The second value tends to
zero as $ k \rar \infty $.
\end{rem}

{\bf Proof.}   It is easy to satisfy to the conditions of the
theorem, since $ A_{s(k)} \rightarrow_k 0 $ and we can decrease
the diameters of the sets $ D_s $ when $ q_{s} = A_{s} B_{s} <
\frac{1}{2} $ and the point $ x_k $ can be considered as a new
starting point.

At each step  $ k $ there is a majorant function
$\psi_{s,k}(\cdot)$
$$
\psi_{s} (t) = B_{s} t ^ 2 -2 t +2 A_{s}.
$$

Since
$$
\| Q_{s,0} {\Theta} '' _ {s } (x_k) \| \leq \psi_{s} '' (x_k),
$$
the step length  $ \triangle_k = \| x_ {k + 1} - x_k \| $ does not
exceed the step length $ t_{k+1} - t_k $ of the Newton's method
for the equation $ \psi_{s} (t) = 0 $ the solution of which we
denote by $t_s$ thus the following can be written: \be \| x_{k+1}
- x_k \| \leq t_{k+1} - t_k . \label{equilibpoints19c} \ee

For the existence of the majorant equation $ \psi_{s} (t) = 0, t
\in \mathbb {R}, $ for the operator equation $ \Theta (x) = 0 $,
as it follows from the Taylor formula  XVII.2.5
\cite{kantorovichakilov}, it is sufficient that the following
integral inequality is correct
$$
\| \int_ {x_k} ^ {x_ {k + 1}} Q_0 {\Theta}''_ {s} (x) (x_ {k + 1}
-x, \cdot) dx \| \leq \int_{t_k}^{t_{k + 1}} c_0 \psi_{s} '' (t)
(t_{k + 1}-t) dt
$$
for  big enough $ k $, which is correct if
$$
\| Q_0 {\Theta}''_ {s} (x) \| \, \| x_{k + 1}-x_k \| \leq c_0
\psi_{s}''(t) (t_{k + 1}-t_k).
$$
Let us denote by
$$
c_{s,k}=-\frac{1}{\psi_{s}'(t_k)}, \,\, A_{s,k}= c_{s,k}
\psi_{s}(t_k),\,\,
$$
$$
B_{s,k}=c_{s,k}\psi''_{s}(t_k)=2 B_{s} c_{s,k}, \,\,
q_{s,k}=B_{s,k}A_{s,k}.
$$
Let us note that \be t_{k+1} - t_k = -
\frac{\psi_{s}(t_k)}{\psi'_{s}(t_k)} = A_{s,k}, \,\,\, k=0,1,
\dots \label{equilibpoints20} \ee According to the Taylor
expansion for a second degree polynomial we have
$$
A_{s,k}=c_{s,k}\psi_{s}(t_k)=c_{s}\psi_{s}(t_{k-1}+A_{s,k-1})=
$$
$$
=c_{s,k} \left [  \frac{1}{2}\psi''_{s}(t_{k-1})A_{s,k-1}^2
+\psi'_{s}(t_{k-1}) A_{s,k-1}+\psi_{s}(t_{k-1})\right ]=
$$
$$
=c_{s,k} \left[
B_{s}A_{s,k-1}^2-\frac{A_{s,k-1}}{c_{s,k-1}}+\frac{A_{s,k-1}}{c_{s,k-1}}\right
]=c_{s,k}B_{s}A_{s,k-1}^2=
$$
$$
=\frac{1}{2} \frac{c_{s,k}}{c_{s,k-1}}2 B_{s}c_{s,k-1} A_{s,k-1}^2
= \frac{1}{2}\frac{c_{s,k}}{c_{s,k-1}}B_{s,k-1}A_{s,k-1}^2
$$
However,
$$
\frac{c_{s,k}}{c_{s,k-1}}=\frac{\psi_s'(t_k)}{\psi_s(t_{k-1})}=
$$
\be
=\frac{\psi_s(t_k-1)+\psi'_s(t_k-1)A_{s,k-1}}{\psi_s(t_{k-1})}=
1-B_{s,k-1}A_{s,k-1}=1-q_{s,k-1}. \label{equilibpoints21} \ee
Therefore, \be A_{s,k}=\frac{1}{2}
\frac{B_{s,k}A_{s,k-1}^2}{1-q_{s,k-1}}
=\frac{A_{s,k-1}}{2}\frac{q_{s,k-1}}{1-q_{s,k-1}}.
\label{equilibpoints22} \ee By analogy, from
(\ref{equilibpoints21}) we obtain
$$
B_{s,k}=2 c_{s,k} B_s=2 B_s c_{s,k-1}
\frac{c_{s,k}}{c_{s,k-1}}=\frac{B_{s,k-1}}{1-q_{s,k-1}}.
$$
From here \be q_{s,k}=B_{s,k}A_{s,k}=\frac{1}{2} \frac{B_{s,k-1}
A_{s,k-1} q_{s,k-1}}{(1-q_{s,k-1})^2}=\frac{1}{2}\left [
\frac{q_{s,k-1}}{1-q_{s,k-1}} \right ]^2. \label{equilibpoints23}
\ee From (\ref{equilibpoints22}) and  (\ref{equilibpoints23}),
taking into account $q_{s,k} \leq \frac{1}{2}$, we obtain the
following estimations \be A_{s,k} \leq q_{s,k-1} A_{s,k}, \,\,\,
q_{s,k} \leq 2 q_{s,k-1}^2 \,\, n=1,2, \dots
\label{equilibpoints24} \ee Consequently, $q_{s,k} \leq
\frac{1}{2} [ 2 q_{s,0} ]^{2^k} = \frac{1}{2} [ 2 q_{s} ]^{2^k}$
$$
A_{s,k} \leq q_{s,k-1} A_{s,k-1} \leq q_{s,k-1} q_{s,k-2}
A_{s,k-2} \leq  \dots q_{s,k-1} q_{s,k-2} \dots  q_{s,0} A_{s,0},
$$
in which $q_{s,0}=q_s$ and $A_{s,0}=A_s$.

From here and (\ref{equilibpoints19c}), (\ref{equilibpoints20}) we
obtain
$$
\| x_{k+1} - x_k \| +\| x_{k+2} - x_{k+1} \| + \dots  \leq
(t_{k+1} - t_{k})+(t_{k+2} - t_{k+1})+\dots \leq
$$
\be \leq \frac{1}{2^k} [2 q_s]^{2^k -1} A_s \leq
 \frac{1}{2^k} [2 q_s]^{2^k} \frac{A_s}{q_s} \leq
\frac{1}{2^k} [2 q]^{2^k} C, \label{equilibpoints25} \ee since
$$
\frac{A_s}{q_s}=\frac{1}{B_s} \leq C
$$
and $B_s $ is an upper bound for the norm of the second
derivatives and can only increase as $s \rar \infty$. Passing to
the limit on $k=k(s) \rar \infty $ in (\ref{equilibpoints25}), we
obtain the inequality (\ref{equilibpoints19b}).

We will use the modified Newton method $\Delta x_k=-
(\tilde{\Theta}_{s,0}'(x_0))^{-1} \tilde{\Theta}_{s,k}(x_k)$ for
solving the equality $\Theta_{s}(x) =0 $. We denote  the obtained
sequence as $ \{ x'_k \} $. Suppose
$$
 \varphi_s(t)=t+c_0 \psi_s(t),
$$
where $c_{0} =  -\frac{1}{\psi'_s(t_0)} = \frac{1}{2}$ and $t_0=0
$. Let us replace the modified Newton method for the equation
$\psi_s(t)=0 $ with the equation
$$
t=\varphi_s(t), \,
$$
and we will solve it by using the successive approximations
method. Suppose
$$
t^*_s = \varphi_s (t^*_s).
$$

We can write
$$
t^*_s-t'_k= \varphi_s(t^*_s)-
\varphi_s(t'_{k-1})=\varphi_s'(\tilde t_k)(t^* -t'_{k-1}),\,\,
\tilde t_k=\frac{t^*+t'_{k-1}}{2}.
$$
However,
$$
\varphi_s'(t)=1+c_o \psi_s'(t)= B_s t,
$$
so that
$$
\varphi'(\tilde t_k)=B_s \tilde t_k \leq B_s t^* = 1-\sqrt{1-2
q_s}.
$$
Therefore,
$$
t^*_s - t'_k \leq [ 1-\sqrt{1-2 q_s}] ( t^*_s - t'_{k-1}).
$$
We can obtain the similar inequality for $ t^*_s - t'_{k-1} $.
Consequently,
$$
t^*_s - t'_k \leq [1-\sqrt{1-2 q_s}]^k (t^*_s - t'_0)=
\frac{A_s}{q_s} [1-\sqrt{1-2 q_s}]^{k+1}.
$$
The inequality
$$
\| x^*_s - x'_k \| \leq  t^*_s - t'_k ,
$$
similar to the inequality (\ref{equilibpoints19c}), is correct for
the modified Newton's method. Using this inequality, we get
$$
\| x^*_s - x'_k \| \leq \frac{A_s}{q_s} [1-\sqrt{1-2 q_s}]^{k+1}.
$$
Passing to the limit in $s$ and considering
$$
\frac{A_s}{q_s}=\frac{1}{B_s} \leq C, \,\, q_s \leq q
<\frac{1}{2},
$$
in which $C$ is a constant for all $s$,  we obtain the inequality
(\ref{equilibpoints19d}).

The theorem is proved.

\vspace{0.5cm}

\section{\bf Conclusion}

\vspace{0.5cm}

We propose a method for finding equilibrium points as the limit
points of a sequence obtained by applying the numerical method
described above. The coordinate descent method slowly converges to
an equilibrium point in the general case, but by changing the
initial points, one can obtain all equilibrium points with minimal
intermediate calculations.

A method for finding J. Nash equilibrium points using the matrices
of second mixed derivatives (generalized matrices of second mixed
derivatives) of the original functions is suggested. Such methods,
under certain conditions, converge much faster than the coordinate
descent method, but require more calculations at each step.

To speed up the convergence of the method, it is proposed to
decrease consistently the diameter of the set $ D_m $ on which the
integration is performed, and the step length of the optimization
process. We give the rules for successive decrease of the diameter
of the set $ D_m $ and the step length. The Kantorovich theorem is
used to estimate the convergence rate.

\end{document}